# Self-Similar Structures of Nontransitive Dice Sets: Examples of Nested "Rock-Paper-Scissors" Relations Based on Numbers from The Lo Shu Magic Square

Alexander Poddiakov

**Abstract.** Nontransitive dice are dice beating one another in a cyclic way: die *A* wins die *B*, *B* wins *C*, and *C* wins *A* (like in a rock-paper-scissors game). In this article, it has been shown that a structure of mutual wins of 3 nontransitive dice (with numbers equal to numbers from The Lo Shu magic square) can be repeated at least twice in different scales:
  - in the nontransitive relations between three sets, each of which consists of three nontransitive dice (total 9 dice); and
  - in the nontransitive relations between three sets, each of which consists of three nontransitive subsets, each of which consists of three nontransitive dice (total 27 dice).
  In other words, structures of nontransitive superiority relations can be self-similar.
  Aims of the future study can be:
  - to show opportunities for building self-similar structures of nontransitive relations of arbitrary depths of nestedness; and
  - to design a recursive algorithm of filling the structures with the appropriate numbers. Perhaps a geometrical presentation of these numbers forms a fractal structure.

**1. INTRODUCTION.** In 1963, Martin Gardner [8] gave an example from Leo Moser. Moser had discovered the following paradox hidden in The Lo Shu, ancient Chinese magic square (Table 1).

| A | 4 | 9 | 2 |
|---|---|---|---|
| B | 3 | 5 | 7 |
| C | 8 | 1 | 6 |

Table 1: Strengths of chess players equal to numbers from The Lo Shu and forming an intransitive cycle of wins of teams *A, B* and *C*.

> "Let row *A* be a team of three chess experts with the playing strengths of 4, 9, and 2 respectively. Rows *B* and *C* are two other teams, with playing strengths as indicated. If teams *A* and *B* play a round-robin tournament, in which every player of one team plays once against every player of the other team, team *B* will win five games and team *A* will win four. Clearly team *B* is stronger than *A*. When team *B* plays team *C*, *C* wins five games and loses four, so that *C* is obviously stronger than *B*. What happens when *C*, the strongest team, plays *A*, the weakest? Work it out yourself! Team *A* is the winner by five to four! Which, then, is the strongest team? The paradox brings out the weakness of round-robin play in deciding the relative strengths of teams" [Ibid., p. 116-117].

Let me use the following illustration with lengths of sticks. Let us consider 3 sets of sticks of different lengths. We hold a round-robin tournament in which each stick from one set is compared with each stick from another set. The winner is the longest stick in this pair. Let the lengths of sticks be equal to numbers from The Lo Shu. Let us consider patterns of wins between the sticks' sets (Figures 1-3).

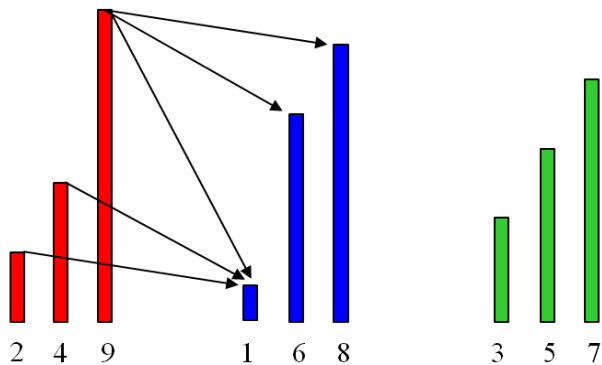

Figure 1: Red sticks "beat" blue ones 5 out of 9 times.

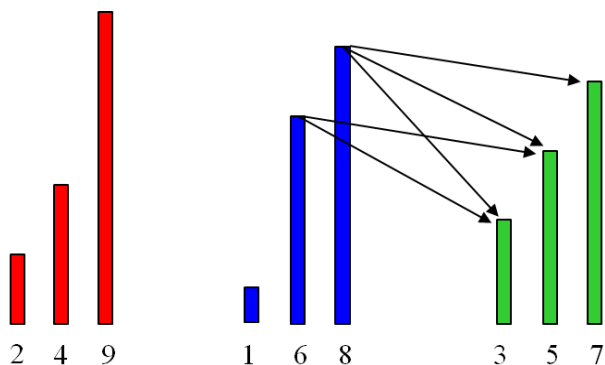

Figure 2: Blue sticks "beat" green ones 5 out of 9 times.

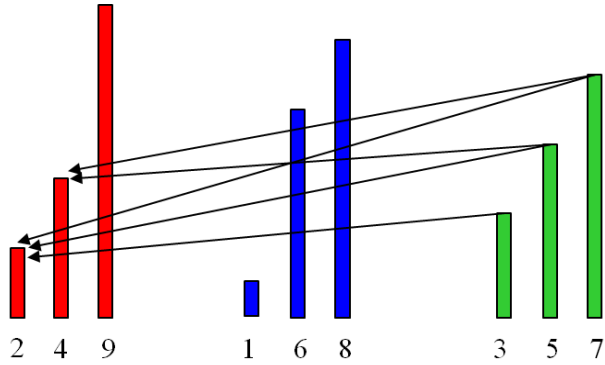

2 4 9      1 6 8      3 5 7

Figure 3: Green sticks "beat" red ones 5 out of 9 times.

One can see that the blue sticks "beat" the green ones 5 out of 9 times, the red sticks "beat" the blue ones 5 out of 9 times, and the green sticks "beat" the red ones 5 out of 9 times.

Teams *A*, *B*, and *C* from Moser's example and the sticks' sets form so called intransitive cycles of dominance (superiority). The intransitive cycle of dominance (superiority) is characterized by such binary relations between A, B, and C that A is superior to B, B is superior to C, and C is superior to A (i.e., A>B>C>A—in contrast with transitive relations A>B>C). Earlier studies on stochastic (not deterministic—like in Moser's example) intransitivity of dominance had been conducted by Trybuła and Steinhaus and published in 1959 and 1961 [22],[23]. Since then, many mathematical works on the intransitivity of dominance have been published, especially on nontransitive dice beating one another like in a rock-paper-scissors game [1]-[6], [9], [20].

A simple example of 6-sided intransitive dice are three dice with double sets of numbers equal to numbers in the Lo Shu Square [21].

Die A has numbers 2, 2, 4, 4, 9, 9 on its faces.
Die B has numbers 1, 1, 6, 6, 8, 8 on its faces.
Die C has numbers 3, 3, 5, 5, 7, 7 on its faces.

These dice win one another with the probability 5/9 in the "rock-paper-scissors" way: A wins B, B wins C, and C wins A (Figure 4).

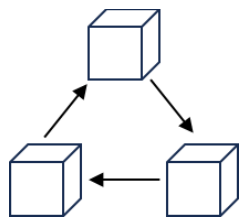

Figure 4: The relationship graph of wins in the nontransitive dice set.

Below I will show that this structure of mutual wins of 3 nontransitive dice can be repeated at least twice in different scales:

(a) in the nontransitive relations between three sets, each of which consists of three nontransitive dice (total 9 dice); and

(b) in the nontransitive relations between three sets, each of which consists of three nontransitive subsets, each of which consists of three nontransitive dice (total 27 dice).

In other words, structures of nontransitive superiority relations can be self-similar. "A self-similar object is one whose component parts resemble the whole. This reiteration of details or patterns occurs at progressively smaller scales and can, in the case of purely abstract entities, continue indefinitely, so that each part of each part, when magnified, will look basically like a fixed part of the whole object. In effect, a self-similar object remains invariant under changes of scale—i.e., it has scaling symmetry" [7].

## 2. NONTRANSITIVE RELATIONS BETWEEN THREE SETS, EACH OF WHICH CONSISTS OF THREE NONTRANSITIVE DICE

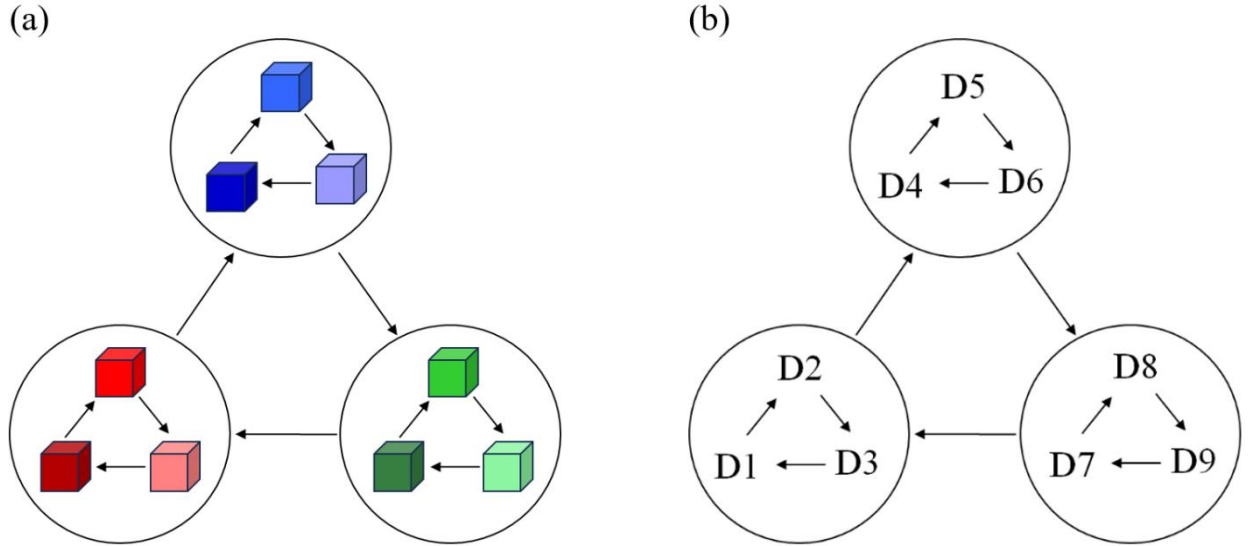

Figure 5: The relationship graph of wins of 3 dice sets, each of which consists of 3 intransitive dice.
(a) Any die of shades of red beats any die of shades of blue, any die of shades of blue beats any die of shades of green, and any die of shades of green beats any die of shades of red. Inside each circle, the darkest die beats the die of intermediate brightness, the die of intermediate brightness beats the lightest die, and the lightest die beats the darkest one.
(b) The graph with numeration of the dice.

Let us consider dice with the following numbers on their faces (Tables 2-4) [18]. (Only half of the numbers on the dice faces are shown in the tables below, and the other numbers are repetitions of the first 3 numbers.)

| The dice of shades of red | Die 1 | 22 | 44 | 99 |
|---|---|---|---|---|
| | Die 2 | 21 | 46 | 98 |
| | Die 3 | 23 | 45 | 97 |

Table 2: Numbers on faces of the dice of shades of red.

| The dice of shades of blue | Die 4 | 12 | 64 | 89 |
|---|---|---|---|---|
| | Die 5 | 11 | 66 | 88 |
| | Die 6 | 13 | 65 | 87 |

Table 3: Numbers on faces of the dice of shades of blue.

| The dice of shades of green | Die 7 | 32 | 54 | 79 |
|---|---|---|---|---|
| | Die 8 | 31 | 56 | 78 |
| | Die 9 | 33 | 55 | 77 |

Table 4: Numbers on faces of the dice of shades of green.

Here any die of shades of red beats any die of shades of blue with a probability 5/9, any die of shades of blue beats any die of shades of green with a probability 5/9, and any die of shades of green beats any die of shades of red with a probability 5/9. In any nested set (a set of dice of shades of a certain color), the darkest die beats the die of intermediate brightness with a probability 5/9, the die of intermediate brightness beats the lightest die with a probability 5/9, and the lightest die beats the darkest one with a probability 5/9.

This nestedness is a new, numerical, example of meta-intransitivity of the 1st level[1]. "Meta-intransitivity is a property of systems, which are in intransitive relationships between one another, and, at the same time, each of them contains its own internal, nested intransitive cycles of superiority" [19]. A number of levels of the nested cycles is the meta-intransitivity level. However, self-similar structures of meta-intransitive relations have not been described yet. (Generally, meta-intransitive cycles can contain nested intransitive cycles which are not self-similar. This text deals with their self-similarity.)

## 3. SELF-SIMILAR STRUCTURES OF 27 DICE AT THE 2ND LEVEL OF META-INTRANSITIVITY

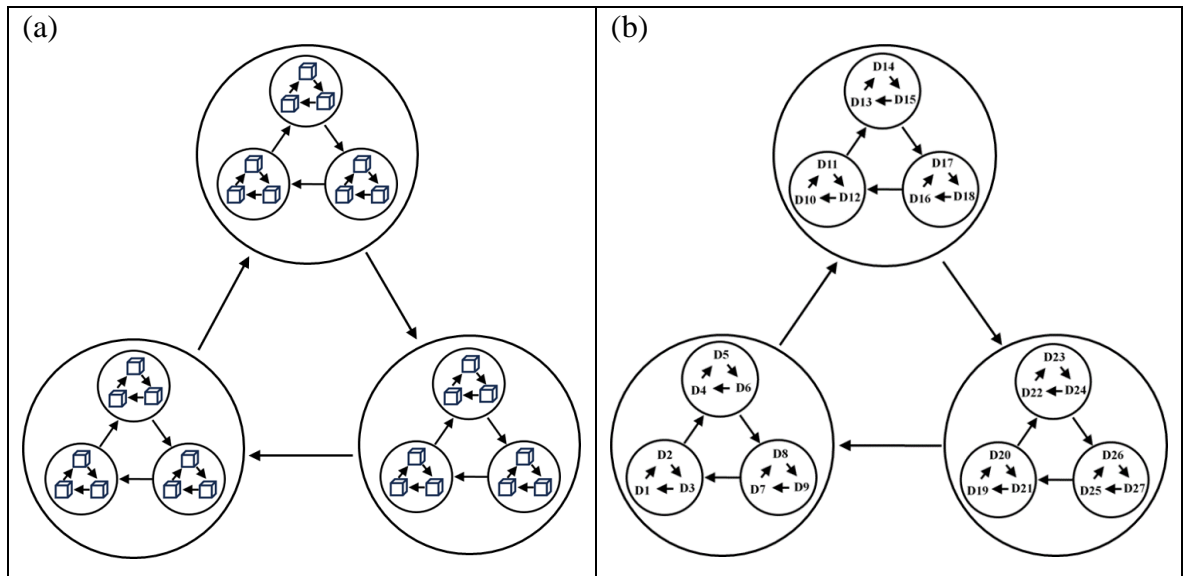

Figure 6: The relationship graph of wins of 3 dice sets, each of which consists of 3 subsets consisting of 3 dice: (a) a general view; (b) the graph with numeration of the dice.

Let us consider 27 dice D1, ..., D27 with the following numbers on their faces (Tables 5-7) [18].

| D1 | 222 | 489 | 954 | | D4 | 299 | 464 | 932 | | D7 | 244 | 412 | 979 |
|---|---|---|---|---|---|---|---|---|---|---|---|---|---|
| D2 | 221 | 488 | 956 | | D5 | 298 | 466 | 931 | | D8 | 246 | 411 | 978 |
| D3 | 223 | 487 | 955 | | D6 | 297 | 465 | 933 | | D9 | 245 | 413 | 977 |

Table 5: Numbers on dice D1, …, D9 from the left large circle in Figure 6.

| D10 | 122 | 689 | 854 | | D13 | 199 | 664 | 832 | | D16 | 144 | 612 | 879 |
|---|---|---|---|---|---|---|---|---|---|---|---|---|---|
| D11 | 121 | 688 | 856 | | D14 | 198 | 666 | 831 | | D17 | 146 | 611 | 878 |
| D12 | 123 | 687 | 855 | | D15 | 197 | 665 | 833 | | D18 | 145 | 613 | 877 |

Table 6: Numbers on dice D10, …, D18 from the top large circle in Figure 6.

| D19 | 322 | 589 | 754 | | D22 | 399 | 564 | 732 | | D25 | 344 | 512 | 779 |
|---|---|---|---|---|---|---|---|---|---|---|---|---|---|
| D20 | 321 | 588 | 756 | | D23 | 398 | 566 | 731 | | D26 | 346 | 511 | 778 |
| D21 | 323 | 587 | 755 | | D24 | 397 | 565 | 733 | | D27 | 345 | 513 | 777 |

Table 7: Numbers on dice D19, …, D27 from the right large circle in Figure 6.

---

[1] In spite of the difference of strict logical notions "nontransitive" and "intransitive", in literature on nontransitive (intransitive) dominance relations, these notions are used as synonyms (e.g., notion "nontransitive dice" means the same as "intransitive dice").

One can see the following.

Any of 9 dice D1, ..., D9 beats any of 9 dice D10, …, D18 with a probability 5/9.
Any of 9 dice D10, …, D18 beats any of 9 dice D19, …, D27 with a probability 5/9.
Any of 9 dice D19, …, D27 beats any of 9 dice D1, ..., D9 with a probability 5/9.

At the same time, each of 3 sets (D1, ..., D9; D10, …, D18; and D19, …, D27) has nested structures of nontransitive relations, like in Section 2. Let us describe them.

In set D1, ..., D9:
any of 3 dice D1, D2, D3 beats any of 3 dice D4, D5, D6 with a probability 5/9,
any of 3 dice D4, D5, D6 beats any of 3 dice D7, D8, D9 with a probability 5/9,
any of 3 dice D7, D8, D9 beats any of 3 dice D1, D2, D3 with a probability 5/9.

In set D10, ..., D18:
any of 3 dice D10, D11, D12 beats any of 3 dice D13, D14, D15 with a probability 5/9,
any of 3 dice D13, D14, D15 beats any of 3 dice D16, D17, D18 with a probability 5/9,
any of 3 dice D16, D17, D18 beats any of 3 dice D10, D11, D12 with a probability 5/9.

In set D19, ..., D27:
any of 3 dice D19, D20, D21 beats any of 3 dice D22, D23, D24 with a probability 5/9,
any of 3 dice D22, D23, D24 beats any of 3 dice D25, D26, D27 with a probability 5/9,
any of 3 dice D25, D26, D27 beats any of 3 dice D19, D20, D21 with a probability 5/9.

At the last level of the nestedness,
for $i = \{1, 4, 7, 10, 13, 17, 20, 23, 25\}$,
Die $i$ beats Die $i+1$ with a probability 5/9,
Die $i+1$ beats Die $i+2$ with a probability 5/9,
Die $i+2$ beats Die $i$ with a probability 5/9.

Thus, we have filled the self-similar structure of nested nontransitive relations presented in Fig. 6 with the appropriate numbers satisfying the conditions of the nestedness at the 1st and 2nd levels of meta-intransitivity. The self-similar patterns of wins in the nontransitive dice sets consisting of 3, 9, and 27 dice are shown in Figure 7.

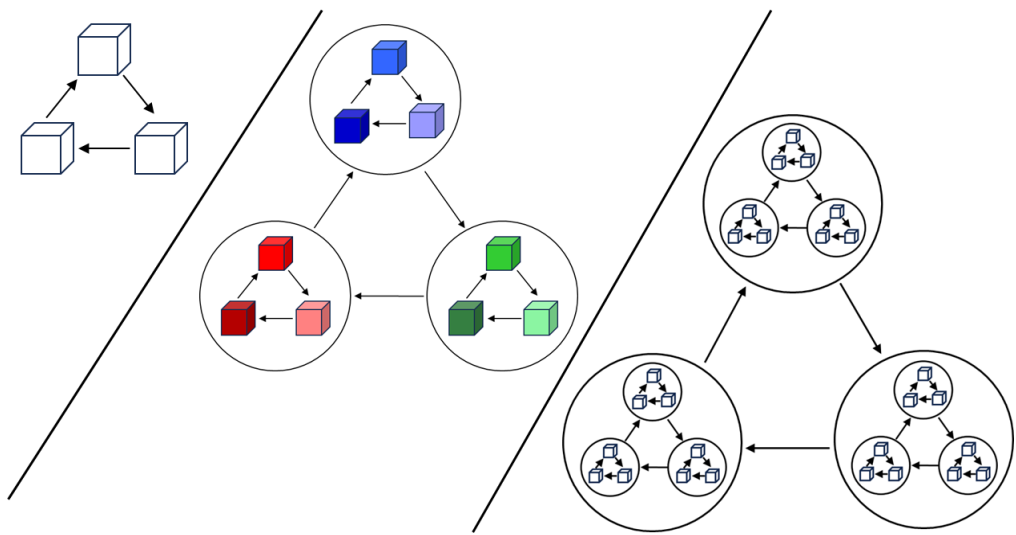

Figure 7: The self-similar patterns of wins in the nontransitive dice sets consisting of 3, 9, and 27 dice.

Aims of the future study can be:

(a) to show opportunities for building self-similar structures of nontransitive relations of arbitrary depths of nestedness; and

(b) to design a recursive algorithm of filling the structures with the appropriate numbers.

Perhaps a geometrical presentation of these numbers forms a fractal structure. A search for laws of fractal structures for sets with different, arbitrarily large initial numbers of nontransitive dice (the case of three nontransitive dice considered above is the simplest one) can be the final point of this study of multiple nested nontransitive relations. (Perhaps numbers normalized to 1 are more suitable for it: 0.211, 0.212, 0.213 etc.)